\documentstyle[11pt,amssymb,amstex]{amsart}
\numberwithin{equation}{section}
\newcommand{\bc}{\begin{center}}
\newcommand{\ec}{\end{center}}
\newcommand{\be}{\begin{equation}}
\newcommand{\ee}{\end{equation}}
\newcommand{\bea}{\begin{eqnarray}}
\newcommand{\eea}{\end{eqnarray}}
\newcommand{\ba}{\begin{array}}
\newcommand{\ea}{\end{array}}

\def\e{{\bf 1}\!\!{\rm I}}

\def\l{\lambda}
\def\O{\Omega}

\def\m{\mu}

\def\o{\omega}
\def\O{\Omega}
\def\nb{\nabla}

\begin{document}

\title[On the "Zero-Two" Law]{On the "Zero-Two" Law for Positive Contractions on
Banach-Kantorovich Lattice $L^p(\nabla,\mu)$}
\author{Inomjon Ganiev}
\address{Inomjon Ganiev\\
Tashkent Railway Institute\\
Tashkent, Uzbekistan} \email{{\tt inam@@tps.uz}}
\author{Farruh Mukhamedov}
\address{Farrukh Mukhamedov\\
Department of Mechanics and Mathematics\\
National University of Uzbekistan\\
Vuzgorodok, 700174, Tashkent, Uzbekistan} \email{{\tt
far75m@@yandex.ru}}

\begin{abstract}
In the present paper we prove the "zero-two" law for positive
contractions of lattices $L^p(\nabla,\mu)$ of Banach-Kantorovich,
constructed by the measure $\mu$ with values in the ring of all
measurable functions.

\end{abstract}

\maketitle

\section{Introduction}

In \cite{W} some properties of the convergence of Banach-valued
martingales was described and there connections with the
geometrical properties of Banach spaces was established too. In
accordance with the development of the theory of
Banach-Kantorovich spaces (see \cite{KVP},\cite{K1},\ \cite{K2},
\cite{G1},\ \cite{G2}) one naturally arises a necessity to study
some ergodic properties of positive contractions and martingales
defined on these Banach-Kantorovich spaces. In \cite{CG} an analog
of individual ergodic theorem for positive contractions of
$L^p(\nb,\m)$ - Banach-Kantorovich space has been established. In
\cite{Ga3} the convergence of martingales on such spaces was
proved.

Let $(X,\Sigma,\m)$ be a measure space and let $L^p(X,\m)$,
$(1\leq p\leq \infty)$ be the usual real $L^p$-space. A linear
operator $T:L^p(X,\m)\to L^p(X,\m)$ is called positive contraction
if for every  $x\in L^p(X,\m)$,$x\geq 0$ we have $Tx\geq 0$ and
$\|T\|_p\leq 1$.

 In \cite{OS} Ornstein and Sucheston  proved that
for any positive contraction $T$ on an $L^1$-space, either
$\|T^n-T^{n+1}\|=2$ for all $n$ or
$\lim\limits_{n\to\infty}\|T^n-T^{n+1}\|=0$. An extension of this
result to positive operators on $L^{\infty}$-spaces was given by
Foguel \cite{F}. In \cite{Z1},\cite{Z2} Zahoropol generalized
these results, called "zero-two" laws, his result can be
formulated as follows:

{\bf Theorem A.} {\it Let $T$ be a positive contraction of
$L^p(X,\m)$, $p>1,p\neq 2$. If for some $m\in
{\mathbb{N}}\cup\{0\}$ the relation
$\bigg\||T^{m+1}-T^m|\bigg\|_p<2$ is valid, then
$$
\lim_{n\to\infty}\|T^{n+1}-T^n\|_p=0.
$$
}

In \cite{KT} this  result was generalized for the K\"{o}the
spaces.

In the present paper we will prove the "zero-two" law for positive
contractions of lattices $L^p(\nabla,\mu)$ of Banach-Kantorovich,
constructed by the measure $\mu$ with values in the ring of all
measurable functions.

\section{Preliminaries}

Let $(\Omega,\Sigma,\lambda)$ be a measurable space with finite
measure, $L_0(\Omega)$ be the algebra of all measurable functions
on $\O$ ( here the functions equal a.e. are identified) and let
$\nb(\O)$ be the Boolean algebra of all idempotents in
$L_0(\Omega)$. By $\nb$ we denote an arbitrary complete Boolean
subalgebra of $\nb(\O)$.

A mapping $\mu : {\nabla}\to L_0(\Omega)$ is called a {\it
$L_0(\Omega)$-valued measure} if
the following conditions are satisfied\\
1) $\mu(e)\geq 0$ for all $e\in{\nabla}$;\\
2) $\mu(e\vee g)=\mu(e)+\mu(g)$, if $e,g\in{\nabla},\ e\wedge g=0$;\\
3) if $e_n\downarrow 0,\ e_n\in{\nabla},\ n=\overline{1,\infty}$,
then $\mu(e_n)\downarrow 0$.

A $L_0(\Omega)$-valued measure is called {\it strictly positive}
if $\allowbreak \mu(e)=\nobreak 0,\ \allowbreak e\in{\nabla}$
implies $e=0$.

In the sequel we will consider strictly positive
$L_0(\Omega)$-valued measure $\m$ with the following property
$\m(ge)=g\m(e)$ for all $e\in\nb$ and $g\in\nb(\O)$.

By $X(\nb)$ we denote an extremal completely non-connected
compact, corresponding to a Boolean algebra $\nb$, and by
$L_0(\nb)$ it is denoted an algebra of all continuous functions on
$X(\nb)$, which take the values $\pm\infty$ on no-where dense sets
in $X(\nb)$ \cite{S}. It is clear that $L_0(\O)$ is a subalgebra
of $L_0(\nb)$.

Following \cite{B},\cite{S} the known scheme of construction of
$L^p$-spaces, it can be defined a space $L^p(\nb,\m)$ as follows
$$
L^p(\nb,\m)=\left\{f\in L_0(\nabla): \int|f|^pd\m - \textrm{exist}
\ \right\}, \ \ \ p\geq 1
$$
here $\m$ is a  $L_0(\O)$-valued measure on $\nb$.

Let $E$ be a linear space over the real field $\mathbb{R}$. By
$\|\cdot\|$ we denote a $L_0(\Omega)$-valued norm on $E$. Then the
pair $(E,\|\cdot\|)$ is called a {\it lattice-normed space (LNS)
over $L_0(\Omega)$}. An LNS $E$ is said to be {\it
$d$-decomposable} if for every $x\in E$ and decomposition
$\|x\|=f+g$ with $f$ and $g$ disjoint positive elements in
$L_0(\Omega)$ there exist $y,z\in E$ such that $x=y+z$ and
$\|y\|=f$, $\|z\|=g$.

Suppose that $(E,\|\cdot\|)$ is an LNS over $L_0(\O)$. A net
$\{x_\alpha\}$ of elements of $E$ is said to be {\it
$(bo)$-converging} to $x\in E$ (in this case we write
$x=(bo)$-$\lim x_\alpha$), if the net $\{\|x_\alpha - x\|\}$
$(o)$-converges to zero in $L_0(\Omega)$ (written as $(o)$-$\lim
\|x_\alpha -x\|=0$). A net $\{x_\alpha\}_{\alpha\in A}$ is called
{\it $(bo)$-fundamental} if $(x_\alpha-x_\beta)_{(\alpha,\beta)\in
A\times A}$ $(bo)$-converges to zero.

An LNS in which every $(bo)$-fundamental net $(bo)$-converges is
called $(bo)$-complete. A {\it Banach-Kantorovich space (BKS) over
$L_0(\Omega)$} is a $(bo)$-complete $d$-decomposable LNS over
$L_0(\Omega)$. It is well known \cite{K1},\cite{K2} that every BKS
$E$ over $L_0(\Omega)$ admits an $L_0(\Omega)$-module structure
such that $\|fx\|=|f|\cdot\|x\|$ for every $x\in E,\ f\in
L_0(\Omega)$, where $|f|$ is the modulus of a function $f\in
L_0(\Omega)$.

It is known \cite{K1} that $L^p(\nb,\m)$ is a BKS over $L_0(\O)$
with respect to $L_0(\O)$-valued norm
$|f|_p=\bigg(\int|f|^pd\m\bigg)^{1/p}$. Moreover, $L^p(\nb,\m)$ is
a module over $L_0(\O)$.

It is natural that these $L^p(\nb,\m)$ spaces should have many of
analogical properties of the classical $L^p$-spaces, constructed
by real valued measures. The proofs of such properties can be
realized by using one of the following methods:

1. Step by step repeating all steps of the known proofs of
classical $L^p$-spaces, accounting specially properties of
$L_0(\O)$-valued measures.

2. Using Boolean-valued analysis, which gives a possibility to
reduce $L_0(\O)$-modulus $L^p(\nb,\m)$ to the classical
$L^p$-spaces, in the corresponding set theory.

3. Representation of $L^p(\nb,\m)$ as a measurable bundle of
classical $L^p$-spaces.

The first methods really is  not effective, since it has to repeat
all known steps of the proofs modifying them to $L_0(\O)$-valued
measures. The second one connected with the use of  drawing an
enough labour-intensive apparatus of Boolean-valued analysis and
its realization requires a huge preparatory work, which connects
with an establishing of intercommunications in ordinary and
Boolean-valued methods for studied objects of the set theory.

More natural way to investigate the properties of $L^p(\nb,\m)$ is
the using of the third method. Since, one has a sufficiently well
explored  theory of measurable decompositions of Banach lattices
\cite{G1}. Therefore, its an effective using gives a possibility
to obtain various properties of BKS\cite{Ga1},\cite{Ga2}.

Let $(\O,\Sigma,\m)$ be the same as above and  let $X$ be a
assisting a real Banach space $X(\o)$ to each point
$\omega\in\Omega$. A {\it section} of $X$ is a function $u$
defined $\l$-almost everywhere in $\O$ that takes values $u(\o)\in
X(\o)$ for all $\o$ in the domain $dom(u)$ of $u$.

Let $L$ be a set of sections. The pair $(X,L)$ is called a {\it
measurable Banach bundle over $\Omega$} if
\begin{enumerate}
   \item[(1)] $\lambda_1 u_1 +\lambda_2 u_2 \in L$ for every $\lambda_1,\lambda_2\in\mathbb{R}$
   and $u_1,u_2\in L$,where\\
 $\lambda_1 u_1 + \lambda_2 u_2 : \omega\in dom(u_1) \cap
dom(u_2) \to \lambda_1 u_1(\omega) + \lambda_2 u_2(\omega)$;
   \item[(2)] the function $\|u\| : \omega\in dom(u) \to
\|u(\omega)\|_{X(\omega)}$ is measurable for every $u\in L$;
    \item[(3)] the set $\{u(\omega) : u\in L, \omega\in dom(u)\}$ is dense in
    $X(\omega)$ for every $\o\in\O$.
\end{enumerate}

A section $s$ is called {\it step-section} if it has the form
$$
s(\omega)=\sum\limits_{i=1}^n \chi_{A_i}(\omega) u_i(\omega),
$$
for some $u_i\in L$, $A_i\in\Sigma$, $A_i\cap A_j=\emptyset$,
$i\neq j$, $i,j=1,\cdots,n$, $n\in\mathbb{N}$, where $\chi_A$ is
the indicator of a set $A$. A section $u$ is called {\it
measurable} if for every $A\in\Sigma$ with $\l(A)<\infty$ there
exists a sequence of step-functions $\{s_n\}$ such that
$s_n(\o)\to u(\o)$ $\l$-almost everywhere on $A$.

Denote by $M(\Omega,X)$ the set all measurable sections, and by
$L_0(\Omega,X)$ the factor space of $M(\Omega,X)$ over the
equivalence relation of equality a.e. Clearly, $L_0(\O,X)$ is an
$L_0(\O)$-module. We denote the equivalence class of an element
$u\in M(\O,X)$ by  $\hat{u}$. The norm of $\hat{u}\in L_0(\O,X)$
is defined as a class of equivalence containing the function
$\|u(\omega)\|_{X(\omega)}$, i.e.
$\|\hat{u}\|=\widehat{(\|u(\o)\|_{X(\o)})}$.

In \cite{G1} it was proved that $L_0(\Omega,X)$ is a BKS over
$L_0(\O)$. Furthermore, for every BKS $E$ over $L_0(\O)$ there
exists a measurable Banach bundle $(X,L)$ over $\O$ such that $E$
is isomorphic to $L_0(\Omega)$.

Put

$${\cal L}^\infty (\Omega,X)=\{u\in M(\Omega,X) :
\|u(\omega)\|_{X(\omega)}\in {\cal L}^\infty(\Omega) \},$$
$$
L^\infty(\Omega,X)=\{\hat{u}\in L_0(\Omega,X) : u\in {\cal
L}^\infty (\Omega,X),\ u\in\hat{u} \},$$ where ${\cal
L}^{\infty}(\O)$ is the set all bounded measurable functions on
$\O$.

In the spaces ${\cal L}^\infty(\Omega,X)$ and $L^\infty(\Omega,X)$
it can be defined real-valued norms as follows $\|u\|_{{\cal
L}^\infty (\Omega,X)} = \sup\limits_\omega
\|u(\omega)\|_{X(\omega)}$ and
$\|\hat{u}\|_{L^\infty(\Omega,X)}=\||\hat{u}|\|_{L^\infty(\Omega)}$,
respectively.

 A BKS $({\cal U},\|\cdot\|)$ is called a {\it
Banach-Kantorovich lattice} if  ${\cal U}$ is a vector lattice and
the norm $\|\cdot\|$ is monotone, i.e.  $|u_1|\leq|u_2|$ implies
$\|u_1\|\leq\|u_2\|$. It is known \cite{K1} that the cone ${\cal
U}_+$ of  positive elements is $(bo)$-closed. Note the space
$L^p(\nabla,\mu)$ is a Banach-Kantorovich lattice \cite{K1}.

Let $X$ be a mapping assisting an $L^p$-space constructed by
real-valued measure $\m_\o$, i.e.  $L^p(\nabla_\omega,\mu_\omega)$
to each point $\o\in\O$ and let
$$
L=\{\sum\limits_{i=1}^n \lambda_i e_i : \lambda_i\in {\mathbb{R}},
\ \ e_i(\omega)\in\nabla_{\omega},\ i=\overline{1,n},\
n\in\mathbb{N}\}$$ be a set of sections. It is known \cite{Ga2}
that the pair $(X,L)$ is a measurable bundle of Banach lattices
and $L_0(\Omega,X)$ is modulo ordered isomorphic to
$L^p(\nabla,\mu)$.

Let $p$ be a lifting in $L^\infty(\Omega)$ (see \cite{G1}). Let as
before $\nb$ be an arbitrary complete Boolean subalgebra of
$\nb(\O)$ and let $\m$ be an $L_0(\O)$-valued measure on $\nb$.
The set of all essentially bounded functions w.r.t. $\m$ from
$L_0(\nb)$ is denoted by  $L^\infty(\nb,\m)$.

It is known \cite{CG} that there exists a mapping $\ell :
L^\infty(\nb,\m)(\subset L^{\infty}(\O,X)) \to {\cal L}^\infty
(\Omega,X)$, which satisfies the following conditions:
\begin{enumerate}
   \item[(1)] $\ell(\hat{u})\in\hat{u}$ for all $\hat{u}$ such that
   $dom(\hat{u})=\O$;
   \item[(2)]
   $\|\ell(\hat{u})\|_{L^p(\nb_{\o},\m_\o)}=p(|\hat{u}|_p)(\omega)$;
   \item[(3)] $\ell(\hat{u}+\hat{v})=\ell(\hat{u})+\ell(\hat{v})$
   for every $\hat{u},\hat{v}\in L^\infty(\nb,\m)$;
   \item[(4)] $\ell(h\cdot\hat{u})=p(h)\ell(\hat{u})$ for every
   $\hat{u}\in L^\infty(\nb,\m),\ h\in L^\infty(\Omega)$;
   \item[(5)] $\ell(\hat{u})(\o)\geq 0$ whenever $\hat{f}\geq 0$;
   \item[(6)] the set $\{\ell(\hat{u})(\omega) : \hat{u}\in
    L^\infty(\nb,\m)\}$ is dense in $X(\omega)$ for all $\omega\in\Omega$;
   \item[(7)]
   $\ell(\hat{u}\vee\hat{v})=\ell(\hat{u})\vee\ell(\hat{v})$
    for every $\hat{u},\hat{v}\in L^\infty(\nb,\m)$.
\end{enumerate}

The mapping $\ell$ is called a  {\it vector-valued lifting} on
$L^\infty(\nb,\m)$ associated with the lifting $p$ (cp.
\cite{G1}).

Let as before $p\geq 1$ and $L^p({\nabla},{\mu})$ be a
Banach-Kantorovich lattice, and let
$L^p(\nabla_\omega,\mu_\omega)$ be the corresponding $L^p$-spaces
constructed by real valued measures. Let $T:L^p({\nabla},{\mu})\to
L^p({\nabla},{\mu})$ be a linear mapping. As usually we will say
that $T$ is {\it positive} if $T\hat{f}\geq 0$ for every
$\hat{f}\geq 0$.

Suppose that $T$ is a  $L_0(\Omega)$-bounded mapping, namely,
there exists a function $k\in L_0(\Omega)$ such that
$|T\hat{f}|_p\leq k |\hat{f}|_p$ for all $\hat{f}\in
L^p(\nabla,\mu)$. In this case we can define an element of
$L_0(\O)$ as follows
$$
\|T\| =\sup\limits_{|\hat{f}|_p\leq\e} |T\hat{f}|_p,
$$
which is called an {\it $L_0(\O)$-valued norm} of $T$. If
$\|T\|\leq\e$ then a mapping $T$ is said to be a {\it
contraction}.

Now we give an example of  nontrivial contraction.

{\bf Example.} Let  $(\Omega,{\nabla},\lambda)$ be a measurable
space with finite measure and let ${\nabla}_0$ be a right Boolean
subalgebra of ${\nabla}$. $\lambda_0$ denotes the restriction of
$\lambda$ onto ${\nabla}_0$. Now let  $E(\cdot | {\nabla}_0)$ be a
conditional expectation from $L_1(\Omega,{\nabla},\lambda)$ onto
$L_1(\Omega,{\nabla}_0,\lambda_0)$. It is clear that
${\mu}(\hat{e})=E(\hat{e} |{\nabla}_0)$ is a strictly positive
$L_1(\Omega,{\nabla}_0,\lambda_0)$-valued measure on ${\nabla}$.

Let ${\nabla}_1$ be another arbitrary  right Boolean subalgebra of
${\nabla}$ such that ${\nabla}_1\supset{\nabla}_0$. By ${\mu}_1$
we denote the restriction of ${\mu}$ onto ${\nabla}_1$. According
to Theorem 4.2.9\cite{K1} there exists an conditional expectation
$T : L_1({\nabla},{\mu})\to L_1({\nabla}_1,{\mu}_1)$ such that it
is positive and maps  $L^p({\nabla},{\mu})$ onto
$L^p({\nabla},{\mu})$ for all $p>1$. Moreover,
$|T\hat{f}|_p\leq|\hat{f}|_p$ for every $\hat{f}\in
L^p({\nabla},{\mu})$ and $T\e=\e$.\\

In the sequel we will need the following

{\bf Theorem 2.1.} {\it Let $T : L^p({\nabla},{\mu})\to
L^p({\nabla},{\mu})$ be a positive linear contraction such that
$T\e\leq\e$. Then for every $\omega\in\Omega$ there exists a
positive contraction $T_\omega : L^p(\nabla_\omega,\mu_\omega)\to
L^p(\nabla_\omega,\mu_\omega)$ such that $T_\omega f(\omega) =
(T\hat{f})(\omega)$ $\l$-a.s. for every $\hat{f}\in
L^p({\nabla},{\mu})$.}

{\bf Proof.} The positivity of $T$ implies that  $|T\hat{f}|\leq T
|\hat{f}|\leq\|\hat{f}\|_\infty \e$ for every $\hat{f}\in
L^\infty({\nabla},{\mu})$, i.e. the operator $T$ maps from
$L^\infty({\nabla},{\mu})$ to $L^\infty({\nabla},{\mu})$ and it is
continuous in norm $\|\cdot\|_\infty$, where
$\|f\|_{\infty}=varisup|f|$. It is easy to see that for
$\hat{f}\in L^\infty({\nabla},{\mu})$ we have $|T\hat{f}|_p\in
L^\infty(\Omega)$ and $|\hat{f}|_p \in L^\infty(\Omega)$. Define a
linear operator $\varphi(\omega)$ from $\{\ell(\hat{f})(\omega) :
\hat{f}\in L^\infty({\nabla},{\mu})\}$ to
$L^p(\nabla_\omega,\mu_\omega)$ as follows
$$
\varphi(\omega)(\ell(\hat{f})(\omega))=\ell(T\hat{f})(\omega), $$
where $\ell$ is the vector-valued lifting in
$L^\infty({\nabla},{\mu})$ associated with the lifting $p$. From
$|T\hat{f}|_p\leq |\hat{f}|_p$ we obtain \bea
\|\ell(T\hat{f})(\omega)\|_{L^p(\nabla_\omega,\mu_\omega)}=
p(|T\hat{f}|_p)(\omega)\leq\nonumber \\
\leq p(|\hat{f}|_p)(\omega) =
\|\ell(\hat{f})(\omega)\|_{L^p(\nabla_\omega,\mu_\omega)}\nonumber
\eea which implies that the operator $\varphi(\omega)$ is correct
defined and it is bounded. Using the fact that
$\{\ell(\hat{f})(\omega) : \hat{f}\in L^\infty({\nabla},{\mu})\}$
is dense in $L^p(\nabla_\omega,\mu_\omega)$ we can extend
$\varphi(\omega)$ to a continuous linear operator on
$L^p(\nabla_\omega,\mu_\omega)$. This extension is denoted by
$T_\omega$.

Show that $T_\omega$ is positive. Indeed, let $f(\omega)\in
L^p(\nabla_\omega,\mu_\omega)$ and $f(\omega)\geq 0$. Then there
exists a sequence $\{\hat{f}_n\}\subset L^\infty({\nabla},{\mu})$
such that $\ell(\hat{f}_n)(\omega)\to f(\omega)$ in norm of
$L^p(\nabla_\omega,\mu_\omega)$. Consider $\hat{g}_n=\hat{f}_n\vee
0$. Then $\hat{g}_n\geq 0$, and  according to the properties of
the vector-valued lifting $\ell$ we infer
$$
\ell(\hat{g}_n)(\omega) = \ell(\hat{f}_n)(\omega) \vee 0 \to
f(\omega)\vee 0 = f(\omega)$$ in norm of
$L^p(\nabla_\omega,\mu_\omega)$. Whence
$$0\leq\ell(T\hat{g}_n)(\omega)=\varphi(\omega)(\ell(\hat{g}_n)(\omega))
\to T_\omega(f(\omega)),$$ this means $T_\omega f(\omega)\geq 0$.
It is clear that $\|T_\omega\|_{p,\o}\leq 1$ and $T_\omega
f(\omega) = (T\hat{f})(\omega)$ a.s. for every  $\hat{f}\in
L^\infty({\nabla},{\mu})$, here $\|\cdot\|_{\infty}$ is the norm
of operator from $L^\infty(\nb_\o,\m_\o)$ to
$L^\infty(\nb_\o,\m_\o)$.

Now let $\hat{f}\in L^p({\nabla},{\mu})$. Since
$L^\infty({\nabla},{\mu})$ is $(bo)$-dense in
$L^p({\nabla},{\mu})$, then there is a sequence
$\{\hat{f}_n\}\subset L^\infty({\nabla},{\mu})$ such that
$|\hat{f}_n - \hat{f}|_p\to ^{(o)}0$. Then $\|f_n(\omega) -
f(\omega)\|_{L^p(\nabla_\omega,\mu_\omega)} \to 0$ for almost all
$\o$. The equality $T\hat{f} = |\cdot|_p\mbox{-}\lim\limits_{n}
T\hat{f}_n$ implies that
$$\|T_\omega f_n(\omega) -
(T\hat{f})(\omega)\|_{L^p (\nabla_\omega,\mu_\omega)} = \|
(T\hat{f}_n) (\omega) - (T\hat{f})(\omega)
\|_{L^p(\nabla_\omega,\mu_\omega)} \to 0 \ \ \textrm{a.s.
$\omega$}
$$
this means that $(T\hat{f})(\omega) = \lim\limits_n T_\omega
f_n(\omega)$ a.s. On the other hand the continuety of $T_\omega$
yields that $\lim\limits_n T_\omega f_n(\omega) = T_\omega
f(\omega)$ a.s. Hence for every $\hat{f}\in L^p({\nabla},{\mu})$
we have $(T\hat{f})(\omega) = T_\omega f(\omega)$ a.s. This
completes the proof.

\section{Main Results}

In this section we will prove an analog of Theorem A formulated in
the introduction.

Before proving main result we give some useful assertions.

 {\bf Proposition 3.1.} {\it Let $T^{(i)} :
L^p({\nabla},{\mu})\to L^p({\nabla},{\mu})$, $i=1,2$ be positive
linear contractions such that $T\e\leq\e$. Then
$$
\|T^{(1)}-T^{(2)}\|(\o)=\|T^{(1)}_\o-T^{(1)}_\o\|_{p,\o}, \ \
\textrm{a.e.}
$$
here as before $\|\cdot\|_{p,\o}$ is the norm of operator from
$L^p(\nb_\o,\m_\o)$ to $L^p(\nb_\o,\m_\o)$. }

{\bf Proof.} According to Theorem 2.1 we have $T^{(i)}_\omega
f(\omega) = (T^{(i)}\hat{f})(\omega)$, $i=1,2$ a.e. for every
$\hat{f}\in L^p(\hat{\nabla},\hat{\mu})$. Using this fact we get
\bea
|(T^{(1)}-T^{(2)})\hat{f}|_p(\o)=\|((T^{(1)}-T^{(2)})\hat{f}(\o)\|_{L^p(\nb_\o,\m_\o)}=\nonumber\\
=\|(T^{(1)}_\o-T^{(2)}_\o)f(\o)\|_{L^p(\nb_\o,\m_\o)}\leq\nonumber\\
\leq \|T^{(1)}_\o-T^{(2)}_\o\|_{p,\o}\|f(\o)\|_{L^p(\nb_\o,\m_\o)}
\nonumber \eea this implies
$$
\|T^{(1)}-T^{(2)}\|(\o)\leq\|T^{(1)}_\o-T^{(2)}_\o\|_{p,\o}, \ \
\textrm{a.e.}\eqno(1)
$$

By similar reasoning we obtain \bea
\|(T^{(1)}_\o-T^{(2)}_\o)f(\o)\|_{L^p(\nb_\o,\m_\o)}=|(T^{(1)}-T^{(2)})\hat{f}|_p(\o)\leq
\nonumber\\
\leq \bigg(\|T^{(1)}-T^{(2)}\||\hat{f}|_p\bigg)(\o)=\nonumber\\
=\|T^{(1)}-T^{(2)}\|(\o)|\hat{f}|_p(\o)=\nonumber\\
=\|T^{(1)}-T^{(2)}\|(\o)\|f_\o\|_{L^p(\nb_\o,\m_\o)},\nonumber
\eea which yields
$$
\|T^{(1)}-T^{(2)}\|(\o)\geq\|T^{(1)}_\o-T^{(2)}_\o\|_{p,\o}. \ \
\textrm{a.e.}
$$

This inequality with (1) implies the required equality. This
competes the proof.

{\bf Proposition 3.2.} {\it Let $T^{(i)} : L^p({\nabla},{\mu})\to
L^p({\nabla},{\mu})$, $i=1,2$ be positive linear contractions such
that $T\e\leq\e$. Then
$$
\bigg\||T^{(1)}_\o-T^{(2)}_\o|\bigg\|_{p,\o}\leq
\bigg\||T^{(1)}-T^{(2)}|\bigg\|(\o), \ \ \textrm{a.e.}
$$
here $|\cdot|$ means the modulus of an operator. }

{\bf Proof.} Using the formula
$$
|Ax|\leq |A||x|,
$$
where $A:E\to E$ is a linear operator and $E$ is a vector lattice
(see \cite{V}, p.231), we have
$$
|(T^{(1)}_\o-T^{(2)}_\o)g(\o)|\leq\bigg(|T^{(1)}-T^{(2)}||\hat{g}|\bigg)(\o) \
\ a.e.
$$
for every $\hat{g}\in L^p(\nb,\m)$.

Let $|\hat{g}|\leq |\hat{f}|$, where $\hat{f}\in L^p(\nb,\m)$,
then $|g(\o)|\leq |f(\o)|$. This implies
$$
|(T^{(1)}_\o-T^{(2)}_\o)g(\o)|\leq\bigg(|T^{(1)}-T^{(2)}||\hat{f}|\bigg)(\o)
\ \ a.e.
$$
Now using the formula
$$
|A|x=\sup_{|y|\leq x}|Ay|,
$$
where $A$ is as above and $x\geq 0$ (see \cite{V}, p.231) we get
$$
|T^{(1)}_\o-T^{(2)}_\o||f(\o)|=\sup_{|g(\o)\leq|f(\o)|}|(T^{(1)}_\o-T^{(2)}_\o)g(\o)|\leq
\bigg(|T^{(1)}-T^{(2)}||\hat{f}|\bigg)(\o).
$$
Then the monotonicity of the norm $\|\cdot\|_{L^p(\nb_\o,\m_\o)}$
implies \bea
\bigg\|\bigg(|T^{(1)}_\o-T^{(2)}_\o||{f}|\bigg)(\o)\bigg\|_{L^p(\nb_\o,\m_\o)}\leq
\left\|\bigg(|T^{(1)}-T^{(2)}||\hat{f}|\bigg)(\o)\right\|_{L^p(\nb_\o,\m_\o)}=\nonumber\\
=\bigg||T^{(1)}-T^{(2)}||\hat{f}|\bigg|_p(\o)\leq\nonumber\\
\leq \left(\bigg\||T^{(1)}-T^{(2)}|\bigg\||\hat{f}|_p\right)(\o)=\nonumber\\
=\bigg\||T^{(1)}-T^{(2)}|\bigg\|(\o)|\hat{f}|_p(\o)=\nonumber\\
=\bigg\||T^{(1)}-T^{(2)}|\bigg\|(\o)\|f(\o)\|_{L^p(\nb_\o,\m_\o)}.\nonumber
\eea

Thus \bea
\bigg\||T^{(1)}_\o-T^{(2)}_\o|\bigg\|_{p,\o}=\sup_{\|f(\o)\|_{L^p(\nb_\o,\m_\o)}\leq
1}\bigg\||T^{(1)}_\o-T^{(2)}_\o||f(\o)|\bigg\|_{L^p(\nb_\o,\m_\o)}\leq\nonumber
\\ \leq\bigg\||T^{(1)}-T^{(2)}|\bigg\|(\o) \ \ a.e.\nonumber
\eea

The next theorem is an analog of  theorem in \cite{Z2} for
positive contractions of $L^1(\nb,\m)$.

{\bf Theorem 3.3.} {\it Let $T: L^1({\nabla},{\mu})\to
L^1({\nabla},{\mu})$, be  a positive linear contraction such that
$T\e\leq\e$. If for some $m\in{\mathbb{N}}\cup\{0\}$ the relation
$\|T^{m+1}-T^m\|<2\e$ is valid. Then
$$
(o)-\lim_{n\to\infty}\|T^{n+1}-T^{n}\|=0.
$$}

{\bf Proof.} According to Theorem 2.1 there exist positive
contractions $T_\o:L^1(\nb_\o,\m_\o)\to L^1(\nb_\o,\m_\o)$ such
that $(T\hat{f})(\o)=T_\o(f(\o))$ a.e. From Proposition 3.1 we get
$\|T^{m+1}_\o-T^m_\o\|_{p,\o}=\|T^{m+1}-T^m\|(\o)$ a.e. The
condition of theorem implies $\|T^{m+1}_\o-T^m_\o\|_{p,\o}<2$ a.e.
Hence the contractions $T_\o$ are satisfy the condition of Theorem
1.1. \cite{OS} a.e., therefore
$$
\lim_{n\to\infty}\|T^{n+1}_\o-T^n_\o\|_{p,\o}=0 \ \ \ a.e.
$$
As $\|T^{n+1}_\o-T^n_\o\|_{p,\o}=\|T^{n+1}-T^n\|(\o)$ a.e. we
obtain that
$$
\lim_{n\to\infty}\|T^{n+1}-T^n\|(\o)=0 \ \ a.e.
$$
therefore
$$
(o)-\lim_{n\to\infty}\|T^{n+1}-T^n\|=0.
$$

The theorem is proved.

Now  we can prove the following theorem, which is an analog of
Theorem A for the Banach-Kantorovich lattice
$L^p({\nabla},{\mu})$.

{\bf Theorem 3.4.} {\it Let $T: L^p({\nabla},{\mu})\to
L^p({\nabla},{\mu})$, $p>1, p\neq 2$ be a positive linear
contraction such that $T\e\leq\e$. If for some
$m\in{\mathbb{N}}\cup\{0\}$ the relation
$\bigg\||T^{m+1}-T^m|\bigg\|<2\e$ is vaild. Then
$$
(o)-\lim_{n\to\infty}\|T^{n+1}-T^{n}\|=0.
$$}

The proof uses the similar argument as Theorem 3.3, here instead
of Proposition 3.1 it should be used Proposition 3.2.\\[2mm]

\end{document}